\documentclass{article}[11]%
\usepackage{latexsym,amssymb,amsmath,amscd,amsthm,amsxtra}
\usepackage{amsmath}
\usepackage{amsfonts}
\usepackage{amssymb}
\usepackage{graphicx}%
\newtheorem{thm}{Theorem}[section]

\newtheorem{dfn}[thm]{Definition}

\newtheorem{rmk}[thm]{Remark}

\newtheorem{prop}[thm]{Property}

\newcommand{\reals}{\mathbb {R}}

\newcommand {\pf}{\noindent{\bf Proof.}\ }

\newcommand{\arrows}{\,\lower1pt\hbox{$\longrightarrow$}\hskip-.24in\raise2pt
             \hbox{$\longrightarrow$}\,}
\title{A Counter Example of Invariant Deformation Quantization}
\author{Xiang Tang}
\begin{document}
\maketitle
\begin{abstract}
In this note, we will show one example of hamiltonian Lie algebra
action which has no invariant star product.
\end{abstract}

\section{Introduction}
Quantization of a hamiltonian system with symmetries is an important
and difficult problem in physics and mathematics. In the deformation
quantization formulation \cite{bffls:deformation}, this problem can
be phrased as follows: given a hamiltonian Lie group action on a
symplectic manifold, does there exist a star product containing(see
Definition \ref{dfn:inv-star}) the information of the group action?

Since the very early time of deformation quantization,
Lichnerowicz has considered this question(see \cite{lich:star} and
references therein). Lichnerowicz in \cite{lich:connection} showed
that if a homogeneous space $G/H$ admits an invariant linear
connection, the symplectic manifold $T^*(G/H)$ admits an invariant
Vey $\star$-product.

In literature, there are various definitions of a star product. We
fix our star product to the following one.

\begin{dfn}
\label{dfn:star} Let $(M, \omega)$ be a symplectic manifold. A
star product on $C^{\infty}(M)$ is an associative product $\star$
on $C^{\infty}(M)[[\hbar]]$ with the following properties:
\begin{enumerate}
\item the coefficients $c_k(x)$ of the product
$$
c(x, \hbar)=a(x, \hbar)\star b(x, \hbar)=\sum_{k=0}^{\infty}\hbar^k
c_k(x)
$$
depend on $a_i, b_j$ and their derivatives $\partial^{\alpha}a_i,
\partial^{\beta}b_j$.
\item the leading term $c_{0}(x)$ is equal to the usual commutative
product of functions $a_0(x)b_0(x)$.
\item the $\star$-product satisfies
$$
[a, b]=a\star b-b\star a=-i\{a_0, b_0\}+\cdots,
$$
where $\{\ ,\ \}$ means the Poisson bracket of functions and dots
mean higher order terms of $\hbar$.
\end{enumerate}
\end{dfn}

In this note, we will write a star product as $f\star
g=\sum_{r=0}^{\infty}\hbar^rC^r(f, g)$, where $C^r$ is a local
bidifferential operator. Next we recall the definition of a Vey
$\star$-product.
\begin{dfn}
\label{dfn:vey-star}Let $\nabla$ be a symplectic connection on $(M,
\omega)$. A $\star$ product is called a Vey $\star^n$-product if the
principal symbol of the differential operator $C^r$ is identical to
$$
P^r_{\nabla}(f, g)=\omega^{i_1j_1}\cdots
\omega^{i_rj_r}\nabla_{i_1}\cdots\nabla_{i_r}f\nabla_{j_1}\cdots\nabla_{j_r}g,\
\ for\ all\ f,\ g\in C^{\infty}(M),
$$
for all $r\leq n$.
\end{dfn}

At the beginning of this section, when describing the question of
quantization with symmetries, we have been very vague by using the
word ``containing". In literatures, there are several related
notions of invariant and covariant star products. In this paper, we
will focus on the following invariant star product from
\cite{ar:invariance}.
\begin{dfn}
\label{dfn:inv-star}For a hamiltonian Lie group $G$ action on a
symplectic manifold $(M, \omega)$, a $\star$ product is called
strongly $G$ invariant\footnote{In short, we will just say ``$G$
invariant" star product in this note.} if:
$$
x\cdot (f\star g)=(x\cdot f)\star (x\cdot g), \ \ for\ all\ x\in
G,\ f,\ g\in C^{\infty}(M).
$$
\end{dfn}

Looking at the infinitesimal Lie algebra $\frak{g}$ action and
$J:\frak{g}\to C^{\infty}(M)$ the dual of the momentum map, we have
$$
\{ J(X), f\star g\}=\{J(X), f\}\star g+f\star\{ J(X), g\},
$$
for all $X$ in $\frak{g}$, $f,\ g$ in $C^{\infty}(M)$.

From Definition \ref{dfn:vey-star} and \ref{dfn:inv-star}, we can
easily see that if a Vey$^2-$product is $G$-invariant, then the
corresponding symplectic connection is also $G-$invariant.
Therefore, Lichnerowicz's result is also necessary for the existence
of an invariant Vey$^2$-product. A $G$-invariant Vey$^2-$product
exists if and only there is an invariant symplectic connection.

In Fedosov's construction \cite{Fe:deformation} of star products on
a symplectic manifold, it is obvious that the existence of an
invariant connection implies the existence of an invariant
symplectic connection and therefore the existence of an invariant
star product. With this and the integration trick, for any
hamiltonian compact Lie group action, we can construct an invariant
connection and therefore an invariant star product.

The existence of invariant star products leads to the study of
quantum momentum map and reduction theory.  Xu in \cite{xu:moment}
introduced and studied the theory of quantum momentum map. In
\cite{fe:s1-quotient} and \cite{fe:quotient}, Fedosov used his
quantization method to study quantum Marsden-Weinstein reduction of
a compact hamiltonian Lie group action. Bordemann, Herbig, and
Waldmann in \cite{bhs:brst} studied BRST cohomology in the framework
of deformation quantization and quantum reduced space.

Recently, in literature, there are many attempts to extend the study
of invariant star products and Xu's quantum momentum map to more
general types of quantization. In \cite{neu:invariant},
M\"uller-Bahns and Neumaier considered star products of wick type;
and in \cite{gutt:invariant}, Gutt and Rawsly investigated natural
star products. All the known results have suggested that the
original idea of Lichnerowicz that the existence of an invariant
star product is closely related to the existence of an invariant
connection is correct.

In the above discussion, we have concentrated on symplectic
manifolds. It is worth mentioning the Poisson version of the
question. The existence of a star product for a general Poisson
manifold was first constructed by Kontsevich(and later Tarmarkin
with a different method) in \cite{ko:deformation} using his
formality theorem. From Kontsevich's original construction, it is
not very obvious to see the conditions needed for the existence of
an invariant star product. Dolgushev in \cite{do:deformation} gave
an alternative construction of the global formality theorem using
Fedosov type resolution and Kontsevich's local formality theorem.
Dolgushev's construction explicitly shows that the existence of an
invariant connection is a sufficient condition for an invariant star
product(also invariant formality theorem). It would be interesting
to look at the Poisson version of quantum momentum maps and BRST
quotients.

It is also worth mentioning that since \cite{cmz:modular} and
\cite{dlo:conformal}, there has been discussion of conformally
invariant symbol calculus and star products. These products are
different from the star product defined in Definition \ref{dfn:star}
that they are highly nonlocal. The study of conformally invariant
quantization is still at its early stage, and we even do not know
whether a conformally invariant quantization always exists. However,
we have seen its interesting relations to other areas of
mathematics. For example, Cohen, Manin, and Zagier in
\cite{cmz:modular} obtained this type of products from considering
deformation of modular forms. In \cite{bty:modular}, we will use
this type of star products(also Fedosov's construction) to
reconstruct Connes and Moscovici's universal deformation formula
\cite{cm:deformation} of the Hopf algebra associated to
codimensional one foliation.

In this note, we will show that there is a hamiltonian Lie algebra
action which has no invariant star product, which can be viewed as
an analog of Van Hove's ``no-go" theorem in invariant deformation
quantization.

In this direction, Arnal, Cortet, Molin, and Pinczon in
\cite{ar:invariance} showed that on some coadjoint orbit
${\mathcal{O}}$ of a nilpotent Lie algebra $\frak{g}$, there is no
$\frak{g}$-invariant $Vey^2-$product by showing that there is no
invariant $\frak{g}-$connection.

What we will do is basically to extend their result to any star
product. Since we are working in full generality, to show that there
is no invariant connection as in \cite{ar:invariance} is not enough
any more. We will study properties of general invariant differential
operators, which will give us enough information to show the
nonexistence of an invariant star product.

\begin{rmk}
\label{rmk:ex}This type of counter examples is at least believed to
exist among experts we have talked to. But we cannot find any
explicit example in literatures. If there is any other examples,
please let us know.
\end{rmk}

\begin{rmk}
\label{rmk:semisimple}On a large class of coadjoint orbits,
invariant star products were constructed in \cite{ak:invariant} and
references therein.
\end{rmk}

\begin{rmk}
\label{rmk:covariant}Weaker than invariant star products, people
have introduced a notion of ``covariant star products"(see
\cite{ar:invariance}). Instead of the keeping the same action, we
allow higher order modification to the group(Lie algebra) action.
The existence and uniqueness of covariant star products are related
to the lower order Lie algebra(Lie group) cohomology(see
\cite{neu:invariant}). This year there are many interesting
activities in this direction. This spring, Kontsevich conjectured
that the automorphism group of the Poisson algebra of polynomial
functions on $R^{2n}$ is naturally isomorphic to the automorphism
group of the corresponding $2n-$dimensional Weyl algebra. And this
summer in IHP, Gorokhovsky, Nest, and Tsygan showed the author a
very interesting construction of their ``stacky star product".
\end{rmk}
{\bf {Acknowledgement}:} The result of this paper was completed
during my Ph. D. study in UC Berkeley. Firstly, I would like to
thank my thesis advisor Alan Weinstein for proposing this question
to me and many helpful comments and suggestions. I also want to
thank Simone Gutt for answering me many questions in emails, and
Gorokhovsky, Nest, and Tsygan for interesting discussion.
\section{Main result}
We look at $(\reals ^2, dx\wedge dy)$ with the Lie algebra $\frak
{g}$ action formed by the Hamiltonian vector fields generated by
\[
x^3,\ x^2,\ x,\ y,\ 1.
\]

$\frak {g}$ is a 5-dim nilpotent Lie algebra\footnote{We can look at
the Lie algebra of the corresponding hamiltonian vector fields,
which has no center.}. By the expression of a star product, we can
easily see that if a $\star$-product is invariant under $\frak {g}$
action, then each $C^r$ of $\star$ has to be $\frak{g}$ invariant,
i.e.
\[
X(C^r (u,v))=C^r(X(u), v)+C^r(u, X(v))\ \ \ \ \ \forall X\in \frak
{g},\ \ u,\ v\in C^{\infty}(M),\ \ r=1,\ 2,\ 3,\ \cdots .
\]
Therefore, in the following, we will first look at properties of
differential operators that are invariant under the $\frak {g}$
action.  Then we will come back to the existence of an invariant
$\star \index{$\star$}$ product.

We write a bidifferential operator as
\[
C_{ij; kl} \index{$C_{ij; kl}$}(\partial _x)^i (\partial _y)^j
\otimes (\partial _x)^k (\partial _y)^l,
\]
where we have used the Einstein summation convention.

\begin{prop}
\label{prop:inv-opt} If a bidifferential operator $C_{ij; kl}
\index{$C_{ij; kl}$}(\partial _x)^i (\partial _y)^j \otimes
(\partial _x)^k (\partial _y)^l$ is invariant under the $\frak
{g}$ action, then $C_{ij;kl}$ satisfies the following relations:
\begin{enumerate}
\item $C_{ij;kl}$ are all constants.
\item if $i> l$ or $j<k$, then $C_{ij;kl} =0$;
\item
\[
\begin{array}{lll}
C_{ij;kl}&=-C_{i+1, j-1; k-1, l+1},\ \ &for\ j\geq1,\ k \geq 1; \\
C_{ij;kl}&=-C_{i-1, j+1; k+1, l-1},\ \ &for\ i\geq1,\ l \geq 1.
\end{array}
\]
\item
\[
\begin{array}{lll}
C_{ij;kl}&=-C_{i+2, j-1; k-2, l+1},\ \ &for\ j\geq1 ,\ k\geq 2;\\
C_{ij;kl}&=-C_{i-2, j+1; k+2, l-1},\ \ &for\ i\geq2 ,\ l\geq 1.
\end{array}
\]
\end{enumerate}
\end{prop}
$\pf$ We work on each generator of $\frak {g}$.
\begin{enumerate}
\item $1\in \frak {g}$. This part is trivial. Because the
hamiltonian vector field of 1 is 0, every bidifferential operator
is invariant under it. \item $x\in \frak {g}$. The hamiltonian
vector field generated by $x$ is $\partial _y$. If $C_{ij; kl}
\index{$C_{ij; kl}$}(\partial _x)^i (\partial _y)^j \otimes
(\partial _x)^k (\partial _y)^l$ is invariant under $\partial_y$,
then
\[
\begin{array}{ll}
\partial_y (C_{ij; kl} \index{$C_{ij; kl}$}(\partial _x)^i (\partial _y)^j \otimes (\partial _x)^k
(\partial _y)^l)&=C_{ij; kl} \index{$C_{ij; kl}$}(\partial _x)^i
(\partial _y)^{j+1}
\otimes (\partial _x)^k(\partial _y)^l\\
&+C_{ij; kl} \index{$C_{ij; kl}$}(\partial _x)^i (\partial _y)^j
\otimes (\partial _x)^k(\partial _y)^{l+1}.
\end{array}
\]
We expand the left hand side of the above equation, and after
cancellations, we have
\[
\partial _y(C_{ij;kl})=0.
\]
\item $y\in \frak {g}$. Similar to the case of $x$, we get
\[
\partial _x(C_{ij;kl})=0.
\]

From the above, we have $\partial _x(C_{ij;kl})=\partial
_y(C_{ij;kl})=0$ on $\reals ^2$. Therefore, $C_{ij;kl}$ is a
constant. \item $x^2\in \frak {g}$. The Hamiltonian vector field
generated by $x^2$ is $2x\partial _y$. The invariance of $C_{ij;
kl} \index{$C_{ij; kl}$}(\partial _x)^i (\partial _y)^j \otimes
(\partial _x)^k (\partial _y)^l$ under $2x\partial _y$ gives
\begin{equation}
\begin{array}{ll}
\label{x^2} 2x\partial _y(C_{ij; kl} \index{$C_{ij; kl}$}(\partial
_x)^i (\partial _y)^j \otimes (\partial _x)^k(\partial
_y)^l)&=C_{ij; kl} \index{$C_{ij; kl}$}(\partial _x)^i (\partial
_y)^j(2x\partial _y) \otimes (\partial
_x)^k(\partial _y)^l+\\
&+C_{ij; kl} \index{$C_{ij; kl}$}(\partial _x)^i (\partial _y)^j
\otimes (\partial _x)^k(\partial _y)^l(2x\partial _y).
\end{array}
\end{equation}

Setting $x=0$ in the above equation, we get
\begin{equation}
\label{eq:simplification} C_{ij;kl}(\partial _x)^{i-1}(\partial
_y)^{j+1}\otimes (\partial _x)^k(\partial _y)^l+C_{ij; kl}
\index{$C_{ij; kl}$}(\partial _x)^i (\partial _y)^j \otimes
(\partial _x)^{k-1} (\partial _y)^{l+1}=0,
\end{equation}
where the first term exists when $i>0$, and the second term exists
when $k>0$.

\begin{enumerate}
\item We look at terms of the form $(\partial _x)^i\otimes (\partial
_x)^{k}{\partial_y}^l$. It is easy to find that the first term of
Equation (\ref{eq:simplification}) does not have this kind of term
since its existence requires $j$ to be greater than or equal to 1.
From this, we have
\[
C_{i0;kl}=0\ \ \ \ \ \forall k>0.
\]
\item Next, we look at terms of the form  $(\partial _x)^i (\partial _y)^j
\otimes (\partial _x)^k$. Arguments like those above show that
\[
C_{ij;k0}=0\ \ \ \ \ \forall i>0.
\]
\item If $j>0, l>0$, by (\ref{x^2}), we get
\[
C_{i+1,j-1;kl}(\partial_x)^i (\partial _y)^j\otimes (\partial
_x)^k(\partial _y)^l +C_{ij;k+1, l-1}(\partial_x)^i (\partial
_y)^j\otimes (\partial _x)^k(\partial _y)^l=0.
\]

This shows that
\[
C_{i+1, j-1; kl}+C_{ij;k+1, l-1}=0 .
\]

Therefore,
\begin{enumerate}

\item if $j>0,\ k>0$,
\[
C_{ij;kl}=-C_{i+1, j-1;k-1, l+1};
\]
\item if $i>0,\ l>0$,
\[
C_{ij;kl}=-C_{i-1, j+1; k+1, l-1}.
\]

According to (a), and iteration using i. of (c), we get that if
$j<k$, $C_{ij;kl}=0$. Similarly, by (b) and ii. of (c), we get
that if $i>l$, then $C_{ij;kl}=0$.
\end{enumerate}
\end{enumerate}
\item $x^3$. The Hamiltonian vector field generated by $x^3$ is $3x^2 \partial _y$.

As in the arguments for $x^2$, we get that
\begin{enumerate}
\item if $k>1$, $C_{i0;kl}=0$; \item if $i>1$, $C_{ij;k0}=0$;
\item if $j\geq 1$ and $l\geq 1$,
\[
C_{i+2, j-1; k, l}+C_{ij; k+2, l-1}=0.
\]

We can rewrite it as the following,
\begin{enumerate}
\item if $j\geq 1$ and $k\geq 2$,
\[
C_{ij;kl}=-C_{i+2, j-1; k-2, l+1};
\]
\item if $i\geq 2$ and $l\geq 1$,
\[
C_{ij;kl}=-C_{i-2, j+1; k+2, l-1}. \ \ \ \Box
\]
\end{enumerate}
\end{enumerate}
\end{enumerate}

With above preparation, we prove the following theorem.
\begin{thm}
For the Hamiltonian $\frak {g}$ action on $(\reals ^2,\ dx\wedge
dy)$, there is no geometrically $\frak {g}$ invariant $\star
\index{$\star$}$ product.
\end{thm}
$\pf$ We prove the theorem by contradiction. Assume that there is
a $\star \index{$\star$}$ product of $(\reals ^2, dx \wedge dy)$
of the form
\[
\sum _{r\geq 0} \hbar ^r C^r ,
\]
which is geometrically $\frak {g}$ invariant.

For each $r>0$, by the assumption of locality, we can write
\[
C^r= C^r _{ij;kl}(\partial _x)^i (\partial _y)^j \otimes (\partial
_x)^k(\partial _y)^l .
\]

According to the associativity of $\star \index{$\star$}$ for the
$\hbar ^2$-term and comparing the corresponding coefficients, we
have that for any $f,\ g,\ h\in C^{\infty}(M)$,
\begin{equation}
\label{eq:assoc-2} C^2(fg,h)+C^1(C^1(f, g), h)+C^2(f, g)h=C^2(f,
gh)+C^1(f, C^1(g,h))+fC^2(g, h).
\end{equation}
\begin{enumerate}
\item We look at the coefficient of the term $f_{yy}g_x h_x$.

\begin{itemize}
  \item On the left hand side of equation (\ref{eq:assoc-2}).
\begin{enumerate}
\item $C^2(fg, h)$. It can possibly contribute the term $C^2
_{12;10}$. But according to the conclusion of Proposition
\ref{prop:inv-opt} that if $i>l$, then $C_{ij;kl}=0$, we have $C^2
_{12;10}=0$. Therefore, $C^2(fg, h)$ has no term of the form
$f_{yy}g_x h_x$.
\item $C^1(C^1(f,g), h)$. There are two $C^1$. As we have $h_x$
term, the outside $C^1$ has to be of the form $C^1_{ij;10}$.
According the result of Proposition \ref{prop:inv-opt}, that if
$i>l$, then $C_{ij;kl}=0$, we have that there are only two
possibilities for the outside $C^1$:
\[
\begin{array}{ll}
C^{1}_{01;10}&and\ \ \ \ \ C^1 _{02;10}.
\end{array}
\]

If the outside $C^1$ contributes $C^1_{01;10}$, then as all the
$C^1 _{ij;kl}$ are constant, the inside one also has to contribute
$C^1_{01:10}$.  Therefore, there is a contribution of
$(C^1_{01;10})^2$.

If the outside $C^1$ has $C^1_{02;10}$, then the inside $C^1$ can
only contribute $C^1 _{00;10}$, but from Proposition
\ref{prop:inv-opt}, it has to be 0, because $j<k$.

So the second term has only one contribution which is
$(C_{01;10})^2$.

\item $C^2 (f,g)h$. Because in this term there is no derivative
respect to $h$, this term can not contribute anything.
\end{enumerate}

In summary, the left hand side of the above equation can only
contribute $(C_{01;10})^2$ to the coefficient of $f_{yy}g_x h_x$.

\item On the right hand side of equation (\ref{eq:assoc-2}).
\begin{enumerate}
\item $C^2(f, gh)$. It can only possibly contribute $C^2
_{02;20}$. But according to Proposition \ref{prop:inv-opt},
\[
C^2 _{02;20}=-C^2_{21;01}.
\]
But from $i>l$, we know $C^{2}_{21;01}=0$. Therefore, there is no
contribution of this term.

\item $C^1(f,C^1(g,h))$. By comparing the number derivatives of $f$, we know that the outside
$C^1$ has to be of the form $C^1 _{02;kl}$. As the differential of
$g$ and $h$ are all respect to $x$, there are three possibility
for the outside $C^1$.
\[
\begin{array}{lll}
C^1_{02;00},&C^1 _{02;10},&and\ \ \ \ C^1 _{02;20}.
\end{array}
\]

In the following, we will show that all three of them do not have
any contribution.
\begin{enumerate}
\item $C^2_{02;00}$. Then the inside $C^1$ has to be of the form
$C^1_{10;10}$. This is $0$ according to Proposition
\ref{prop:inv-opt}. \item $C^2_{02;10}$. Then the inside $C^1$ has
to be of the form $C^1_{10;00}$ or $C^1_{00;10}$, which are both
$0$ because of Proposition \ref{prop:inv-opt}. \item
$C^2_{02;20}$. From the previous calculation, we know that $C^2
_{02;20}=0$.
\end{enumerate}

\item $fC^2(g,h)$. Because this term has no derivative of f, there is no contribution of this term.
\end{enumerate}

In all, total in both sides of equation (\ref{eq:assoc-2}), there
is only one contribution of the term $f_{yy}g_x h_x$, which is
$(C^1 _{01;01}) ^2$. Therefore, we have
\[
C^1 _{01;10}=0 .
\]

\end{itemize}
\item We look at the coefficient of $f_{xx} g_y h_y$.

\begin{itemize}
\item On the left hand side of equation (\ref{eq:assoc-2}).
\begin{enumerate}
\item $C^2 (fg,h)$. The only possible contribution is $C^2
_{21;01}$. But according to Proposition \ref{prop:inv-opt}, $C^2
_{21;01}=0$. \item $C^1(C^1(f, g),h)$. By comparing the
derivatives of $h$, we get that the  outside $C^1$ has to be of
the form $C^1 _{ij;01}$. As $i$ has to be less than or equal to
$1$, otherwise this term is $0$ according to proposition
\ref{prop:inv-opt},  we know that there are four possibilities;
\[
\begin{array}{llll}
C^1 _{10;01},& C^1 _{11;01},&C^1 _{01;01},&and\ \ \ \ C^1
_{00;01}.
\end{array}
\]
In the following, we will show that except for $C^1 _{10;01}$, the
other three cases have no contributions.
\begin{enumerate}
\item $C^1 _{10;01}$. In this case, the inside $C^1 $ also has to
be of the form $C^1 _{10;01}$. The contribution of this term is
$(C^1 _{10;01})^2$. \item $ C^1 _{11;01}$. Then the inside $C^1$
has to be of the form $C^1 _{10;00}$, but this has to be $0$
because $i>l$. So this term has no contribution. \item $C^1
_{01;01}$. Then the inside $C^1$ has to be of the form
$C^1_{20;00}$. This also has to be $0$, because $i>l$. This term
again has no contribution. \item $C^1 _{00;01}$. Then the inside
$C^1$ has to be of the form $C^1 _{20;01}$. This is $0$ for the
same reason as the $C^1 _{20;00}$.
\end{enumerate}
\item $C^{2}(f,g)h$. This has no contribution, because there is no
derivative on $h$.
\end{enumerate}

\item On the right hand side of the relation.

\begin{enumerate}
\item $C^2(f, gh)$. The only possible contribution of $C^2(f, gh)$
is of the form  $C^2_{20;02}$. This has to be $0$, because $ C^2
_{20;02}=C^2 _{01;21}=0$. \item $C^{1}(f, C^1(g, h))$. Comparing
the part of $f$, we know that the outside $C^1$ has to be of the
form $C^{1}_{20,kl}$. As $i$ has to be less than or equal to $l$,
the outside $C^2 $ has to be of the form $C^2 _{20;02}$, which is
0.
\end{enumerate}
In conclusion, total in both sides of equation (\ref{eq:assoc-2}),
there is only one contribution $(C ^1 _{10;01})^2$ for term
$f_{xx}g_yh_y$. Therefore $C^1 _{10;01}=0$.
\end{itemize}
\end{enumerate}

We have shown that $C^{1}_{10;01}$ and $C^1 _{01;10}$ are both
$0$. But on the other hand, from
\[
[u,v]=u\star \index{$\star$} v-v\star \index{$\star$} u=-i\hbar \{
u, v\}+o(\hbar),
\]
we have \[C^1_{10;01}-C^1 _{01;10}=-i.\]

If $C^1 _{10;01}=C^1 _{01;10}=0$, the above  equality can not be
true. So we get a contradiction.

Therefore, there is no geometrically $\frak {g}$ invariant star
product on $(\reals^2, dx\wedge dy)$. $\Box$

Xiang Tang\\
Department of Mathematics\\
University of California, Davis\\
One shields Ave., Davis, CA\\
(xtang@math.ucdavis.edu)

\begin{thebibliography}{99}
\bibitem{ak:invariant}
Alekseev, A., Lachowska, A., Invariant *-products on coadjoint
orbits and the Shapovalov pairing, {\em arxiv:math.QA/0308100}.

\bibitem{ar:invariance}
Arnal, D., Cortet, C., Molin, P. and Pinczon, G., Covariance and
geometrical invariance in *quantization, {\em J. Math. Phys.}24
(2), 276-283.

\bibitem{bty:modular}
Bielivasky, P., Tang, X., and Yao, Y., in preparation.

\bibitem{bffls:deformation}
Bayen, F., Flato, M., Fronsdal, C., Lichnerovicz, A. and
Sternheimer, D. Deformation theory and quantization, {\em Ann.
Physics} {\bf 111}(1978), 61-151.

\bibitem{bhs:brst}
Bordemann, M., Herbig, H., and Waldmann, S. {\em BRST cohomology
and phase space reduction in deformation quantization}, Comm.
Math. Phys.  210  (2000),  no. 1, 107--144.

\bibitem{borde:covariant}
Bordemann, M. (Bi)modules, morphismes et r¨¦duction des
star-produits: le cas symplectique, feuilletages et obstructions,
{\em Arxiv:math.QA/0403334}.

\bibitem{cmz:modular}
Cohen, P., Manin, Y., and Zagier, D., Automorphic
pseudodifferential operators,  {\em Algebraic aspects of
integrable systems,  17--47, Progr. Nonlinear Differential
Equations Appl.}, 26, Birk\"user Boston, Boston, MA, 1997.

\bibitem{cm:deformation}
Connes, A., Moscovici, H., Rankin-Cohen brackets and the Hopf
algebra of transverse geometry,  {\em Mosc. Math. J.}  4  (2004),
no. 1, 111--130, 311.

\bibitem{do:deformation}
Dolgushev, V., Covariant and Equivariant Formality Theorems, to
appear in {\em Adv. Math., Vol. 191}, 1 (2005) 147-177.

\bibitem{dlo:conformal}
Duval, C., Lecomte, P., and Ovsienko, V., Conformally equivariant
quantization: existence and uniqueness.  {\em Ann. Inst. Fourier
(Grenoble)}  49  (1999),  no. 6, 1999--2029.

\bibitem{fe:s1-quotient}
Fedosov, B., Reduction and eigenstates in deformation quantization,
{\em Pseudo-differential calculus and mathematical physics,
277--297, Math. Top.}, 5, Akademie Verlag, Berlin, 1994.

\bibitem{Fe:deformation}
Fedosov, B., {\em Deformation Quantization and Index Theory},
Akademie Verlag, 1996.

\bibitem{fe:quotient}
Fedosov, B., Non-abelian reduction in deformation quantization. {\em
Lett. Math. Phys.} 43 (1998), no. 2, 137--154.

\bibitem{gutt:invariant}
Gutt, S. and Rawnsley, J., Natural star products on symplectic
manifolds and quantum moment maps,  {\em Lett. Math. Phys.}, 66
(2003), no. 1-2, 123--139.

\bibitem {ko:deformation}
Kontsevich, M., Deformation quantization of Poisson manifolds, I,
{\em  Lett. Math. Phys.}. 00:1-61, 2004.

\bibitem{lich:connection}
Lichnerowicz, A., Connexions et existence de $\star $-produits sur
une vari\'et\'e symplectique, (French. English summary) {\em C. R.
Acad. Sci. Paris S\'er}, A-B 291 (1980), no. 7, A463--A467.

\bibitem{lich:star}
Lichnerowicz, A., D\'eformations d'alg\'ebres associ\'ees \'aune
vari\'et\'e symplectique (les $*\sb{\nu }$-produits), (French),
{\em Ann. Inst. Fourier (Grenoble) }, 32 (1982), no. 1, xi--xii,
157--209.

\bibitem{neu:invariant}
M$\ddot{u}$ller-Bahns, M. and Neumaier, N., Some remarks on $
g$-invariant Fedosov star products and quantum momentum mappings,
{\em J. Geom. Phys. 50 (2004)},  no. 1-4, 257--272.

\bibitem{xu:moment}
Xu, P., Fedosov $*$-products and quantum momentum maps, {\em Comm.
Math. Phys.}, 197 (1998), no. 1, 167--197.
\end{thebibliography}
\end{document}